\documentclass[11pt]{article}
\usepackage{amsmath, amssymb, amsthm, calc}

\title{Klein Polyhedra and Lattices with Positive Norm Minima.}
\author{Oleg N. German \thanks{ This research was supported by
                                RFBR (grant $\textup N^\circ$ 06--01--00518),
                                INTAS (grant $\textup N^\circ$ 03--51--5070) and
                                grant of the President of Russian Federation
                                $\textup N^\circ$ MK--6370.2006.1.}}
\date{}

\theoremstyle{definition}
\newtheorem{definition}{Definition}

\theoremstyle{remark}
\newtheorem{remark}{Remark}

\theoremstyle{plain}
\newtheorem{theorem}{Theorem}
\newtheorem{lemma}{Lemma}

\newtheorem{statement}{Statement}

\newtheorem*{corollary_no_counter}{Corollary}
\newtheorem*{littlewood}{Littlewood conjecture}
\newtheorem*{oppenheim}{Oppenheim conjecture on linear forms}
\newtheorem*{sailed_oppenheim}{Reformulated Oppenheim conjecture}

\renewcommand{\vec}[1]{\mathbf{#1}}
\renewcommand{\geq}{\geqslant}
\renewcommand{\leq}{\leqslant}
\renewcommand{\dim}{\textup{dim\,}}
\newcommand{\aff}{\textup{aff}}

\renewcommand{\dim}{\textup{dim\,}}
\newcommand{\conv}{\textup{conv}}
\newcommand{\ext}{\textup{ext\,}}
\newcommand{\interior}{\textup{int\,}}
\newcommand{\starv}{\textup{St}_{\vec v}}
\newcommand{\sail}{\Pi}

\newcommand{\La}{\Lambda}

\newcommand{\R}{\mathbb{R}}
\newcommand{\Z}{\mathbb{Z}}

\newcommand{\N}{\mathbb{N}}

\newcommand{\pilog}{\pi_{\textup{log}}}
\newcommand{\partition}{\mathfrak P}
\newcommand{\faces}{\mathfrak B}

\begin{document}

  \maketitle

  \begin{abstract}
    A Klein polyhedron is defined as the convex hull of nonzero lattice points inside an
    orthant of $\R^n$. It generalizes the concept of continued fraction. In this paper facets
    and edge stars of vertices of a Klein polyhedron are considered as multidimensional
    analogs of partial quotients and quantitative characteristics of these ``partial
    quotients'', so called determinants, are defined. It is proved that the facets of all the
    $2^n$ Klein polyhedra generated by a lattice $\La$ have uniformly bounded determinants if
    and only if the facets and the edge stars of the vertices of the Klein polyhedron
    generated by $\La$ and related to the positive orthant have uniformly bounded
    determinants.
  \end{abstract}

  \section{Introduction}

  In this paper we give a complete proof of the results announced in
  \cite{german_norm_minima_II}. Here we investigate one of the most natural multidimensional
  geometric generalizations of continued fractions, the so--called \emph{Klein polyhedra}.

  Continued fractions admit a rather elegant geometric interpretation
  (see \cite{erdos_gruber_hammer}), which can be described as follows. Given a number
  $\alpha$: $0<\alpha<1$, consider a two--dimensional lattice $\La_\alpha$ with basis vectors
  $(1,1-\alpha)$ and $(0,1)$. The convex hull of the nonzero points of the lattice $\La_\alpha$
  with nonnegative coordinates is called a \emph{Klein polygon}. The integer lengths of the
  Klein polygon's bounded edges are equal to the respective partial quotients of the
  number $\alpha$ with odd indices, and the integer angles between pairs of adjacent edges
  are equal to the partial quotients with even indices. The \emph{integer length} of a
  segment with endpoints in $\La_\alpha$ is defined as the number of lattice points contained
  in the interior of this segment plus one. And the \emph{integer angle} between two such
  segments with a common endpoint is defined as the area of the parallelogram spanned by them
  divided by the product of their integer lengths, or in other words, the index of the
  sublattice spanned by the primitive vectors of $\La_\alpha$ parallel to these two segments.

  If an arbitrary two--dimensional lattice is considered, then there obviously appear two
  numbers with their partial quotients describing the combinatorial structure of the
  corresponding Klein polygon.

  The multidimensional generalization of this construction was proposed more than a century
  ago by F.\,Klein (see \cite{klein}). Let $\La\subset\R^n$ be an $n$--dimensional lattice
  with determinant $1$.

  \begin{definition}
    The convex hulls of the nonzero points of the lattice $\La$ contained in each orthant are
    called \emph{Klein polyhedra} of the lattice $\La$.
  \end{definition}

  In this paper, we consider only irrational lattices $\La$, i.e. we assume that the
  coordinate planes contain no lattice points except the origin $\vec 0$. Then, as shown
  in \cite{moussafir_A_polyhedra}, a Klein polyhedron $K$ is a generalized polyhedron, which
  means that its intersection with an arbitrary bounded polyhedron is itself a polyhedron.
  Hence the boundary of $K$ is in this case an $(n-1)$--dimensional polyhedral surface
  homeomorphic to $\R^{n-1}$, consisting of convex $(n-1)$--dimensional (generalized)
  polyhedra, with each point in it belonging only to a finite number of these polyhedra. Some
  of the faces of $K$ can be unbounded, but only if the lattice, dual to $\La$, is not
  irrational (see \cite{german_norm_minima_I}).

  \begin{definition}
    The boundary $\sail$ of a Klein polyhedron $K$ is called a \emph{sail}.
  \end{definition}

  \begin{definition}
    Let $F$ be a face of $K$ of dimension $k$. We call $F$

    a) a \emph{vertex} of $K$, if $k=0$,

    b) an \emph{edge} of $K$, if $k=1$,

    c) a \emph{facet} of $K$, if $k=n-1$.
  \end{definition}

  A few years ago Vladimir Arnold posed a question (see \cite{arnold_mccme},
  \cite{arnold_preface}) which local affine invariants of a sail are sufficient to
  reconstruct the lattice. This question in its initial formulation remains unanswered.
  However in the current paper we establish a connection between some local invariants of a
  sail and the property of a lattice to have positive norm minimum.

  In the two--dimensional case two neighboring Klein polygons have very much in common, for
  the integer lengths of edges of one of them equal the integer angles between the
  correspondent edges of another one (see, for instance, \cite{korkina_2dim}). Due to this
  fact many statements concerning continued fractions admit ``dual'' formulations: we can use
  only integer lengths of edges, and then we will have to consider all the four Klein
  polygons, or we can use both integer lengths of edges and integer angles between adjacent
  edges, and then we may content ourselves with only one Klein polygon. The main result of
  this paper (Theorem\,\ref{t:the_main_one}) gives an example of a statement on Klein
  polyhedra in an arbitrary dimension admitting such a ``dual'' formulation.

  \section{Formulation of the main result}

  In this section the main result is formulated, which is a
  multidimensional generalization of a well--known statement that a number is badly
  approximable if and only if its partial quotients are bounded. Recall that a number
  $\alpha$ is called \emph{badly approximable} if there exists a constant $c>0$ such that,
  for all integer $p$ and natural $q$, the following inequality holds:
  \[ |q\alpha-p|\geq\frac{c}{q}\ . \]
  In terms of Klein polygons of the lattice $\La_\alpha$ from the beginning of the first
  section, this means exactly that the area
  $ \{ \, \vec x\in\R^2\ |\ x_2>0 \text{ and } |x_1x_2|<c\, \} $
  does not contain any point of $\La_\alpha$.

  Thus, it is natural to consider the property of a lattice $\La$ to have a positive
  norm minimum as a multidimensional generalization of the property of a number to be badly
  approximable.

  \begin{definition}
    \emph{The norm minimum} of a lattice $\La$ is defined as
    \[ N(\La)=\inf_{\vec x\in\La\backslash\{\vec 0\}}|\varphi(\vec x)|, \]
    where $\varphi(\vec x)=x_1\ldots x_n$.
  \end{definition}

  We will also need a multidimensional analog of partial quotients. In view of the
  correspondence between partial quotients and integer lengths and angles mentioned in the
  previous section, it is rather natural in the $n$--dimensional case to expect the
  $(n-1)$--dimensional faces of a sail (we will call them \emph{facets}) and the edge stars
  of a sail's vertices to play the role of partial quotients. As a numerical characteristic
  of these multidimensional ``partial quotients'' we will consider their ``determinants''.

  \begin{definition} \label{d:detF}
    Let $F$ be an arbitrary facet of a sail $\sail$ and let $\vec v_1,\ldots,\vec v_m$ be
    the vertices of $F$. Then, we define the \emph{determinant of the facet} $F$ as
    \[ \det F=
       \sum_{1\leq i_1<\ldots<i_n\leq m}|\det(\vec v_{i_1},\ldots,\vec v_{i_n})|\,. \]
  \end{definition}

  \begin{definition} \label{d:detStv}
    Suppose a vertex $\vec v$ of a sail $\sail$ is incident to $m$ edges. Let
    $\vec r_1,\ldots,\vec r_m$ denote the primitive vectors of the lattice $\La$ generating
    these edges. Then, the \emph{determinant of the edge star} $\starv$ of the vertex
    $\vec v$ is defined as
    \[ \det\starv=
       \sum_{1\leq i_1<\ldots<i_n\leq m}|\det(\vec r_{i_1},\ldots,\vec r_{i_n})|\,. \]
  \end{definition}

  It is clear that when $n=2$, i.e. when the sail is one--dimensional,
  the determinants of the sail's edges are equal to the integer lengths
  of these edges, and the determinants of the edge stars of vertices are
  equal to the integer angles between the correspondent edges.

  Note that we can give an equivalent definition of determinants of facets and edge stars in
  terms of Minkowski sum and mixed volume. Recall (see \cite{bonnesen_fenchel},
  \cite{grunbaum}, \cite{mcmullen_shephard}, \cite{ewald}) that the \emph{Minkowski sum} of
  segments $[\vec 0,\vec x_1],\ldots,[\vec 0,\vec x_m]$ (we will need only this most simple
  case) is the set
  \[ \{ \lambda_1\vec x_1+\cdots+\lambda_m\vec x_m\ |\ 0\leq\lambda_i\leq 1 \} \]
  and its (Euclidean) volume is called the \emph{mixed volume} of the segments
  $[\vec 0,\vec x_1],\ldots,[\vec 0,\vec x_m]$. The following simple statement immediately
  gives us the equivalent way of defining determinants of facets and edge stars:

  \begin{statement}
    For any $\vec x_1,\ldots,\vec x_m\in\R^n$ the mixed volume of the segments
    $[\vec 0,\vec x_1],\ldots,[\vec 0,\vec x_m]$ is equal to
    \[ \sum_{1\leq i_1<\ldots<i_n\leq m}|\det(\vec x_{i_1},\ldots,\vec x_{i_n})|\,. \]
  \end{statement}

  Now that all the needed definitions are given we can formulate the main result of this
  paper. It is a part of the following

  \begin{theorem} \label{t:the_main_one}
    Given an irrational $n$--dimensional lattice $\La\subset\mathbb R^n$, the following
    conditions are equivalent:

    \begin{enumerate}
      \item $N(\La)>0$.

      \item The facets of all the $2^n$ sails generated by $\La$ have uniformly bounded
      determinants (i.e. bounded by a constant not depending on a face).

      \item The facets and the edge stars of the vertices of the sail generated by $\La$ and
      related to the positive orthant have uniformly bounded determinants (i.e. bounded by a
      constant not depending on a face or an edge star).
    \end{enumerate}
  \end{theorem}

  The equivalence of $(1)$ and $(2)$ was established in \cite{german_norm_minima_I}. In the
  current paper are proved the two implications $(1)\,\&\,(2)\implies(3)$ and
  $(3)\implies (2)$.

  \begin{remark}
    Actually, what is proved in this paper is a bit stronger than what is formulated in
    Theorem\,\ref{t:the_main_one}. Namely, it is shown that if $N(\La)=\mu>0$ then there
    exists a constant $D$ depending only on $n$ and $\mu$, such that the facets and the edge
    stars of the vertices of the sail $\sail$ have determinants bounded by $D$. And
    conversely, if the facets and the edge stars of the vertices of the sail have
    determinants bounded by a constant $D$ then there exists a constant $\mu$ depending only
    on $n$ and $D$, such that $N(\La)\geq\mu>0$.
  \end{remark}

  \begin{remark}
    To be precise, the definition of a facet's determinant in \cite{german_norm_minima_I} is
    somewhat different from Definition\,\ref{d:detF}. However it is absolutely clear that the
    uniform boundedness of determinants from \cite{german_norm_minima_I} is equivalent to the
    uniform boundedness of determinants from Definition\,\ref{d:detF}.
  \end{remark}

  \section{A relation to the Littlewood and Oppenheim conjectures}

  The following two conjectures are classical:

  \begin{littlewood}
    If $\alpha,\beta\in\R$, then $\inf_{m\in\N}m\|m\alpha\|\|m\beta\|=0$, where
    $\|\cdot\|$ denotes the distance to the nearest integer.
  \end{littlewood}

  \begin{oppenheim}
    If $n\geq 3\ $ and $\ \La\subset\R^n$ is an $n$--dimensional lattice with $N(\La)>0$,
    then $\La$ is algebraic (i.e. similar modulo the action of the group of diagonal
    $n\times n$ matrices
    to the lattice of a complete module of a totally real algebraic field of degree $n$).
  \end{oppenheim}

  Note that the converse of the latter statement is an obvious corollary of the Dirichlet
  theorem on algebraic units (see \cite{borevich_shafarevich}).

  As is known (see \cite{cassels_swinnerton_dyer}), the three--dimensional Oppenheim
  conjecture implies the Littlewood conjecture. In \cite{skubenko_3dim} and
  \cite{skubenko_ndim} an attempt was made to prove the Oppenheim conjecture, however, there
  was an essential gap in the proof. Thus, both conjectures remain unproved.

  Theorem\,\ref{t:the_main_one} allows to reformulate the Oppenheim conjecture as follows:

  \begin{sailed_oppenheim}
    If $n\geq 3$ and a sail $\sail$ generated by an $n$--di\-men\-sional lattice
    $\La\subset\R^n$ is such that all its facets and edge stars of vertices have uniformly
    bounded determinants, then $\La$ is algebraic.
  \end{sailed_oppenheim}

  It follows from the Dirichlet theorem on algebraic units that a sail generated by an
  $n$--di\-men\-sional algebraic lattice has periodic combinatorial structure. The group of
  ``periods'' is isomorphic to $\Z^{n-1}$ and the fundamental domain is bounded. Thus the
  Oppenheim conjecture yields the following corollary: if a sail's facets and edge stars of
  its vertices have uniformly bounded determinants then this sail has a periodic
  combinatorial structure.

  \section{Dual lattices and polar polyhedra}

  In this section we generalize some of the facts from the theory of polar polyhedra to the
  case of Klein polyhedra. We also prove some statements connecting Klein polyhedra of dual
  lattices. It is worth mentioning in this context the book \cite{lachaud_unseen}, where
  similar questions are considered.

  \begin{definition} \label{d:polar_polyhedron}
    Let $P$ be an arbitrary (generalized) $n$--dimensional polyhedron in
    $\R^n,\ \vec 0\notin P$. Then, the \emph{polar polyhedron} for the $P$ is the set
    \[ P^\circ=\big\{\, \vec x\in\R^n\ \big|\ \forall\vec y\in P\ \
       \langle\vec x,\vec y\rangle\geq 1 \,\big\}. \]
  \end{definition}

  The set $P^\circ$ is obviously closed and convex. Hence we can talk about its \emph{faces},
  defined as the intersections of $P^\circ$ with its supporting hyperplanes. We will denote
  by $\faces(P^\circ)$ and $\faces(P)$ the sets of all proper faces of $P^\circ$ and $P$
  respectively.

  Usually the ``inverse'' concept of polarity is considered, i.e. for polytopes containing
  $\vec 0$ in their interior and with the inverse inequality. And it is well known that
  between the boundary complexes of polar polytopes there exists an inclusion--reversing
  bijection (see, for instance, \cite{grunbaum}, \cite{mcmullen_shephard}, \cite{ewald}). We
  are going to need a similar statement concerning Klein polyhedra generated by irrational
  lattices:

  \begin{statement} \label{st:polar_correspondence_Klein_polyhedra}
    Let $\La$ be an irrational $n$--dimensional lattice in $\R^n$ and $K$ be the Klein
    polyhedron generated by $\La$ and related to the positive orthant. Suppose in addition
    that all the faces of $K$ are bounded. \\
    (a) $K^\circ$ is an $n$--dimensional generalized polyhedron with bounded faces. \\
    (b) If $F$ is a proper face of $K$ then the set $F_K^\circ$ defined as
        \[ F_K^\circ=K^\circ\cap\big\{ \vec x\in\R^n\, \big|\, \langle\vec x,\vec y\rangle=1
           \text{ for all\,\ }\vec y\in F \big\} \]
        is a face of $K^\circ$ and
        \[ \dim F_K^\circ=n-1-\dim F. \]
    (c) The mapping
        \[ \begin{array}{rc}
             \beta_K: & \faces(K)\rightarrow\faces(K^\circ) \\
             \beta_K: & F\mapsto F_K^\circ
           \end{array} \]
        is an inclusion--reversing bijection.
  \end{statement}

  To prove Statement\,\ref{st:polar_correspondence_Klein_polyhedra} we will need three
  auxiliary statements. The first one can be proved simply by literal translation of already
  known arguments for polytopes (see \cite{grunbaum}, \cite{mcmullen_shephard}, \cite{ewald})
  to our case, so we leave it without proof:

  \begin{lemma} \label{l:polar_correspondence_literal_translation}
    Let $P$ be an $n$--dimensional polyhedron in $\R^n,\ \vec 0\notin P$, and let
    $\lambda P\subset P$ for all $\lambda\geq 1$. Suppose that $P$ contains no lines. Let
    $\faces'(P)$ and $\faces'(P^\circ)$ denote the sets of all proper faces of $P$ and $P^\circ$
    respectively, whose affine hulls do not contain $\vec 0$.
    \\
    (a) $P^\circ$ is an $n$--dimensional polyhedron, $\vec 0\notin P,\ \lambda P\subset P$
        for all $\lambda\geq 1$ and $P$ contains no lines. \\
    (b) If $F\in\faces'(P)$ then the set $F_P^\circ$ defined as
        \[ F_P^\circ=P^\circ\cap\big\{ \vec x\in\R^n\, \big|\, \langle\vec x,\vec y\rangle=1
           \text{ for all\,\ }\vec y\in F \big\} \]
        is a face of $P^\circ$ and
        \[ \dim F_P^\circ=n-1-\dim F. \]
    (c) The mapping
        \[ \begin{array}{rc}
             \beta_P: & \faces'(P)\rightarrow\faces'(P^\circ) \\
             \beta_P: & F\mapsto F_P^\circ
           \end{array} \]
        is an inclusion--reversing bijection.
  \end{lemma}

  We will also need the following notation: for each $\vec v\in\R^n$ we will denote by
  $H_{\vec v}^+$ and $H_{\vec v}^-$ the half--spaces
  $\{ \vec x\in\R^n\, |\, \langle\vec x,\vec v\rangle\geq 1 \}$ and
  $\{ \vec x\in\R^n\, |\, \langle\vec x,\vec v\rangle\leq 1 \}$ respectively.

  \begin{lemma} \label{l:P_is_a_convex_hull_of_its_vertices}
    If $P$ is an arbitrary generalized polyhedron containing no lines, and all its edges are
    bounded, then it coincides with the convex hull of its vertices.
  \end{lemma}

  \begin{proof}
    It is enough to show that $P$ is contained in the convex hull of its vertices. Since $P$
    contains no lines, there exists $\vec u\in\R^n$ such that the set
    $P_{\vec u}=P\cap H_{\vec u}^-$ is nonempty and compact. Since $P$ is a generalized
    polyhedron, $P_{\vec u}$ is a bounded polyhedron and, hence, coincides with the convex
    hull of its vertices. But all the vertices of $P_{\vec u}$ are either vertices of $P$ or
    lie on edges of $P$, which are bounded. Therefore, $P_{\vec u}$ is contained in the convex
    hull of vertices of $P$. Hence, so is $P$.
  \end{proof}

  \begin{lemma} \label{l:K_circ_equals_K_prime}
    Let $\La$ be an irrational $n$--dimensional lattice in $\R^n$ and $K$ be the Klein
    polyhedron generated by $\La$ and related to the positive orthant. Then, $K^\circ=K'$,
    where
    \[ K'=\bigcap_{\vec v}H_{\vec v}^+ \]
    and the intersection is taken over all the vertices of $K$.
  \end{lemma}

  \begin{proof}
    Since $\La$ is irrational, every edge of $K$ is bounded. Hence by
    Lemma\,\ref{l:P_is_a_convex_hull_of_its_vertices}, $K$ coincides with the convex hull
    of its vertices. The inclusion $K'\subseteq K^\circ$ easily follows from this fact and
    the Carath\'eodory theorem (see \cite{grunbaum}, \cite{danzer_grunbaum_klee}). The
    inclusion $K^\circ\subseteq K'$ is obvious.
  \end{proof}

  Further by $\conv(M)$ we will denote the convex hull of a set $M$.

  \begin{proof}[Proof of Statement\,\ref{st:polar_correspondence_Klein_polyhedra}]
    For each $\vec u\in\R^n$ let us denote by $V_{\vec u}$ the set of vertices $\vec v$ of
    $K$ such that the open interval $(\vec v,\vec u)$ does not have common points with $K$.
    The set $V_{\vec 0}$ obviously coincides with the set of all vertices of $K$. On the
    other hand, the set $V_{\vec u}$ is finite whenever all $u_i$ are strictly positive,
    since we suppose that $K$ does not have unbounded faces.

    Let us consider an arbitrary point $\vec u\notin K$ with strictly positive coordinates
    and denote
    \[ K_{\vec u}=\bigcup_{\lambda\geq 1}\lambda\,\conv(V_{\vec u}). \]
    Since $V_{\vec u}$ is finite,
    \[ K_{\vec u}^\circ=\bigcap_{\vec v\in V_{\vec u}}H_{\vec v}^+. \]
    At the same time for every $\vec w\in V_{\vec 0}\backslash V_{\vec u}$ the interval
    $(\vec w,\vec u)$ has at least one common point with $\conv(V_{\vec u})$, hence, there
    exist $\lambda_{\vec v}\geq 0$ such that
    \[ \sum_{\vec v\in V_{\vec u}}\lambda_{\vec v}>1\ \ \ \ \text{and}\ \ \ \
       \vec w=\vec u+\sum_{\vec v\in V_{\vec u}}\lambda_{\vec v}(\vec v-\vec u). \]
    Therefore $ K_{\vec u}^\circ\cap H_{\vec u}^-\subset H_{\vec w}^+$. Together with
    Lemma\,\ref{l:K_circ_equals_K_prime} this implies that
    \[ K^\circ\cap H_{\vec u}^-=
       \bigcap_{\vec v\in V_{\vec 0}}H_{\vec v}^+\cap H_{\vec u}^-=
       \bigcap_{\vec v\in V_{\vec 0}\backslash V_{\vec u}}H_{\vec v}^+\cap
       K_{\vec u}^\circ\cap H_{\vec u}^-=
       K_{\vec u}^\circ\cap H_{\vec u}^-. \]

    Thus for each $\vec u\notin K$ with strictly positive coordinates the set
    $K^\circ\cap H_{\vec u}^-$ is a polyhedron. This shows that $K^\circ$ is a generalized
    polyhedron. Consider now an arbitrary facet $F$ of $K^\circ$ and a point $\vec u\notin K$
    with strictly positive coordinates such that the facet $F$ has a nonempty intersection
    with the interior of the half--space $H_{\vec u}^-$. As we have just shown,
    $K^\circ\cap H_{\vec u}^-=K_{\vec u}^\circ\cap H_{\vec u}^-$, so the affine hull of $F$
    coincides with the affine hull of some facet of $K_{\vec u}^\circ$. But $K_{\vec u}$
    satisfies the conditions of Lemma\,\ref{l:polar_correspondence_literal_translation},
    hence, there exists a vertex $\vec v$ of $K$ such that the affine hull of $F$ is given by
    the equation $\langle\vec v,\vec x\rangle=1$. Since the lattice $\La$ is irrational, all
    the coordinates of $\vec v$ are strictly positive, and therefore $F$ is compact. This
    proves (a).

    To prove (b) let us consider an arbitrary proper face $F$ of $K$ and a point
    $\vec u\notin K$ with strictly positive coordinates such that the set $V_{\vec u}$
    contains all the vertices of $F$ and such that the set $F_K^\circ$ is contained in
    $H_{\vec u}^-$. Such points exist since $F$ is compact and is contained in the interior
    of the positive orthant. Then, due to the equality
    $K^\circ\cap H_{\vec u}^-=K_{\vec u}^\circ\cap H_{\vec u}^-$ we have
    \[ \begin{array}{rcl}
       F_K^\circ & = & K^\circ\cap\big\{ \vec x\in\R^n\, \big|\,
       \forall\vec y\in F\,\ \langle\vec x,\vec y\rangle=1 \big\}\ \ =
       \\ & &
       K_{\vec u}^\circ\cap\big\{ \vec x\in\R^n\, \big|\,
       \forall\vec y\in F\,\ \langle\vec x,\vec y\rangle=1 \big\}\ \ =\ \
       F_{K_{\vec u}}^\circ\,.
    \end{array} \]
    Applying to $K_{\vec u}^\circ$ Lemma\,\ref{l:polar_correspondence_literal_translation} we
    get that $F_K^\circ$ is a face of $K_{\vec u}^\circ$ and $\dim F_K^\circ=n-1-\dim F$. But
    $K^\circ\subset K_{\vec u}^\circ$, so $F_K^\circ$ is also a face of $K^\circ$, which
    proves (b).

    It remains to show that $\beta_K$ maps $\faces(K)$ onto $\faces(K^\circ)$. Consider an
    arbitrary $F\in\faces(K^\circ)$ and a point $\vec u\notin K$ with strictly positive
    coordinates such that $F$ is contained in the interior of $H_{\vec u}^-$. The existence
    of such points follows from (a). Then, $F$ is also an element of
    $\faces'(K_{\vec u}^\circ)$ and by
    Lemma\,\ref{l:polar_correspondence_literal_translation} coincides with
    $G_{K_{\vec u}}^\circ$ for some $G\in\faces'(K_{\vec u})$. But the affine hull of $G$
    does not contain $\vec u$, hence, $G\in\faces(K)$, and the equality
    $G_K^\circ=G_{K_{\vec u}}^\circ=F$ completes the proof.
  \end{proof}

  From Statement\,\ref{st:polar_correspondence_Klein_polyhedra} and
  Lemma\,\ref{l:P_is_a_convex_hull_of_its_vertices} we get the following

  \begin{corollary_no_counter} \label{cor:K_circ_is_the_conv_of_its_vertices}
    $K^\circ$ coincides with the convex hull of its vertices.
  \end{corollary_no_counter}

  \begin{definition} \label{d:dual_lattice}
    If vectors $\vec x_1,\ldots,\vec x_n$ form a basis of a lattice $\La\subset\R^n$, then
    the lattice $\La^\ast$ with basis $\vec x_1^\ast,\ldots,\vec x_n^\ast$, such that
    $\langle\vec x_i,\vec x_j^\ast\rangle=\delta_{ij}$, is called \emph{dual} for the lattice
    $\La$.
  \end{definition}

  The lattice $\La^\ast$ also generates a Klein polyhedron in the positive orthant. We will
  denote it by $K^\ast$. From the fact that $\langle\vec x,\vec y\rangle\in\Z$ for all
  $\vec x\in\La$ and $\vec y\in\La^\ast$ one can easily deduce the following

  \begin{statement} \label{st:Kast_is_a_subset_of_Kcirc}
    If the boundary of the positive orthant contains no nonzero points of lattices $\La$ and
    $\La^\ast$, then $K^\ast\subseteq K^\circ$.
  \end{statement}

  Note that in case $n=2$ we can write in Statement\,\ref{st:Kast_is_a_subset_of_Kcirc} that
  $K^\ast=K^\circ$. The reason why for $n\geq 3$ the equality should be substituted by an
  inclusion is that the integer distances from facets of $K$ to $\vec 0$ can be greater than
  1. Here the \emph{integer distance} from a facet $F$ to the origin $\vec 0$ is defined as
  \[ \min_{\vec x_1,\ldots,\vec x_n}|\det(\vec x_1,\ldots,\vec x_n)|, \]
  where the minimum is taken over all the $n$--tuples of linearly independent lattice points
  lying in the affine hull of $F$. The following lemma is obvious:
  \begin{lemma} \label{l:if_distance_is_D_then_the_dual_lattice_gets_divided_by_D}
    If the integer distance from a facet $F$ of $K$ to the origin $\vec 0$ equals $D$, then
    the vertex of $K^\circ$, corresponding to $F$, is a point of the lattice
    $(D^{-1})\La^\ast$.
  \end{lemma}

  It is also worth mentioning that when $n=2$ it is actually the equality $K^\ast=K^\circ$
  that implies that the integer lengths of edges of Klein polygons equal the integer angles
  between the correspondent pairs of adjacent edges. In case $n\geq 3$ there is no such good
  relation between edge stars of $K$ and facets of $K^\ast$ (or facets of $K$ and edge stars
  of $K^\ast$). The reason is that $K^\ast$ and $K^\circ$ do not generally coincide. But even
  if we consider a vertex $\vec v$ of $K$ and the facet $F_{\vec v}$ of $K^\circ$
  corresponding to $\vec v$, it is not clear yet how to connect $\det\starv$ and
  $\det F_{\vec v}$. However we will not need an explicit formula connecting these two
  values, the inequality yielded by the following lemma will satisfy our needs.

  \begin{lemma} \label{l:boundedness_of_detFv}
    Let $F_{\vec v}$ be a facet of the polyhedron $K^\circ$ corresponding to a vertex
    $\vec v$ of the polyhedron $K$. Then, $ \det F_{\vec v}\leq(\det\starv)^{n-1} $.
  \end{lemma}

  Before proving Lemma\,\ref{l:boundedness_of_detFv} let us first prove two auxiliary
  statements.

  \begin{lemma} \label{l:boundedness_of_detFv_simplicial_case}
  Let $\vec r_1,\ldots,\vec r_n$ form a basis of $\R^n$ and let $\vec v\in\R^n$ have positive
  coordinates in this basis. For each $i=1,\ldots,n$ let $F_i$ denote
  the simplex $ \conv\big(\{ \vec v,\vec v+\vec r_1,\ldots,\vec v+\vec r_n \}
  \big\backslash\{ \vec v+\vec r_i \}\big) $ and let $\vec w_i$ be the vector such that
  $\langle\vec w_i,\vec x\rangle=1$ for all $\vec x\in F_i$.

  Then,
  \[ |\det(\vec w_1,\ldots,\vec w_n)|=
     \frac{|\det(\vec r_1,\ldots,\vec r_n)|^{n-1}}{\det F_1\ldots\det F_n}\,. \]
  \end{lemma}

  \begin{proof}
    Let $\vec r_1^\ast,\ldots,\vec r_n^\ast$ denote the basis of $R^n$, dual to the basis
    $\vec r_1,\ldots,\vec r_n$, i.e. the vectors such that
    $\langle\vec r_i,\vec r_j^\ast\rangle=\delta_{ij}$. Then,
    $|\vec w_i|\det F_i=|\vec r_i^\ast||\det(\vec r_1,\ldots,\vec r_n)|$, which implies
    that
    \[ |\det(\vec w_1,\ldots,\vec w_n)|=|\det(\vec r_1^\ast,\ldots,\vec r_n^\ast)|
       \dfrac{|\vec w_1|}{|\vec r_1^\ast|}\cdots\dfrac{|\vec w_n|}{|\vec r_n^\ast|}=
       \dfrac{|\det(\vec r_1,\ldots,\vec r_n)|^{n-1}}{\det F_1\ldots\det F_n}\,. \]
  \end{proof}

  We will denote by $\interior P$ and $\ext P$ the relative interior and the vertex set of a
  polyhedron $P$.
  If $M\subset\R^n$ is a finite set and to each point $\vec x\in M$ a positive mass
  $\nu_{\vec x}$ is assigned, then for each subset $M'\subseteq M$ of cardinality
  $\sharp(M')$ we will denote by $c(M')$ the point
  $(\sum_{\vec x\in M'}\nu_{\vec x}\vec x)/\sharp(M')$, i.e. the center of mass of the set
  $M'$.

  \begin{lemma} \label{l:on_triangulations_and_coverings}
    Let $P$ be a convex $(n-1)$--dimensional polyhedron with arbitrary positive masses
    assigned to its vertices. Let $\mathfrak T$ be an arbitrary partition of the relative
    boundary of $P$ into (closed) simplexes with vertices from $\ext P$. Then,
    \[ \interior P=\bigcup_{\Delta\in\mathfrak T}
       \interior(\conv(\Delta\cup\{ c(\ext P\backslash\ext\Delta) \})). \]
  \end{lemma}

  \begin{proof}
    Let $\vec x\in\interior P$. Then there obviously exists a simplex $\Delta\in\mathfrak T$,
    such that $\vec x\in\conv(\Delta\cup\{c(P)\})$. It remains to notice
    that the set $\conv(\Delta\cup\{c(P)\})\cap\interior P$ is contained in
    $\interior(\conv(\Delta\cup\{ c(\ext P\backslash\ext\Delta) \}))$.
  \end{proof}

  \begin{proof}[Proof of Lemma\,\ref{l:boundedness_of_detFv}]
    The action of the group of diagonal $n\times n$ matrices with determinant 1 obviously
    preserves the combinatorial structure of sails equipped with determinants of facets and
    edge stars. Hence we may consider the point $\vec v$ to have equal coordinates
    $v_1=\ldots=v_n$. Suppose $\vec v$ is incident to $m$ edges of the sail $\sail$. Let
    $\vec r_1,\ldots,\vec r_m$ be the primitive vectors of the lattice $\La$ generating these
    edges. Let us consider arbitrary positive numbers $k_1,\ldots,k_m$ such that the vectors
    $\vec r'_i=k_i\vec r_i$ belong to the same hyperplane and denote
    $P=\conv(\vec r'_1,\ldots,\vec r'_m)$. Consider also an arbitrary point
    $\lambda\vec v\in\interior P$.

    Assign masses $k_i^{-1}$ to the points $r'_i$. Then, by
    Lemma\,\ref{l:on_triangulations_and_coverings}, we can renumber the vectors
    $\vec r_1,\ldots,\vec r_m$ (renumbering accordingly the numbers $k_1,\ldots,k_m$ and the
    vectors $\vec r'_1,\ldots,\vec r'_m$) so that
    $\lambda\vec v=\lambda'_0\vec r'_0+\cdots+\lambda'_{n-1}\vec r'_{n-1}$, where the
    $\lambda'_i$ are strictly positive and $\vec r'_0=(\vec r_n+\cdots+\vec r_m)/(m-n+1)$. We
    set $\vec r_0=\vec r'_0(m-n+1),\ \lambda_0=\lambda'_0/(m-n+1)$ and
    $\lambda_i=k_i\lambda'_i$ for $i=1,\ldots,n-1$ and we get that
    $\lambda\vec v=\lambda_0\vec r_0+\cdots+\lambda_{n-1}\vec r_{n-1}$ with strictly positive
    $\lambda_i$.

    Clearly, $F_{\vec v}$ is contained in the $(n-1)$--dimensional polyhedron $F_{\vec r}$
    defined as
    \[ \begin{array}{rl}
         F_{\vec r}=\big\{ \vec x\in\R^n \ \ \big| &
         \langle\vec x,\vec v\rangle=1 \text{ and }
         \langle\vec x,\vec v+\displaystyle\sum_{i=0}^{n-1}\varkappa_i\vec r_i\rangle\geq 1 \\
         & \text{for all } \varkappa_0\geq 0,\ldots,\varkappa_{n-1}\geq 0 \big\}.
       \end{array} \]
    Since the vectors $\vec r_0,\ldots,\vec r_{n-1}$ are linearly independent and all of the
    coefficients $\lambda_i$ are positive, by
    Lemma\,\ref{l:boundedness_of_detFv_simplicial_case},
    \[ \det F_{\vec r}=\frac{|\det(\vec r_0,\ldots,\vec r_{n-1})|^{n-1}}
       {\prod_{i=0}^{n-1}|\det(\vec v,\{\vec r_j\}_{j=0}^{n-1}\backslash\{\vec r_i\})|}\,. \]
    All the factors in the denominator are nonzero integers, so, applying the inclusion
    $F_{\vec v}\subseteq F_{\vec r}$, we obtain the required estimate.
  \end{proof}

  \section{Uniform boundedness of determinants of a sail's facets}

  In this section are given some facts concerning the sails that enjoy the property that the
  determinants of their facets are uniformly bounded (by a constant depending only on sail).
  We will make use of them in the subsequent sections.

  As before, we denote $\varphi(\vec x)=x_1\ldots x_n$. We also denote by $S(F)$ the
  intersection of the affine hull of a sail's facet $F$ and the positive orthant. We will
  need a value characterizing the volume of the convex hull of $S(F)$ and the origin
  $\vec 0$. It is convenient for this purpose to consider the natural extension of
  Definition\,\ref{d:detF} (given only for facets of a sail) to the case of arbitrary convex
  $(n-1)$--dimensional polyhedra and consider the value $\det S(F)$, which in this case is
  obviously equal to the volume of $\conv(S(F)\cup\{\vec 0\})$ multiplied by $n!$.

  In \cite{german_norm_minima_I} the following is proved:

  \begin{theorem}   \label{t:boundedness_of_F_and_SF}
    Suppose that the boundary of the positive orthant contains no points of a lattice $\La$
    except the origin $\vec 0$. Suppose also that the determinants of all the facets of the
    sail $\sail$ generated by $\La$ and related to the positive orthant are bounded by a
    constant $D$. Then there exists a constant $D'$ depending only on $n$ and $D$ such that \\
    (a) $\det S(F)\leq D'$ for each facet $F$ of the sail $\sail$; \\
    (b) $\varphi(\vec v)\geq (D')^{-1}$ for each vertex $\vec v$ of $K^\ast$.
  \end{theorem}

  \begin{lemma}   \label{l:phi_is_bounded_on_sail}
    If the determinants of all the facets of a sail $\sail$ are bounded by a constant $D$,
    then there exists a constant $D'$ depending only on $n$ and $D$, such that
    $\varphi(\vec x)<D'$ for each point $\vec x\in\sail$.
  \end{lemma}

  \begin{proof}
    It is enough to consider a facet $F$ of the sail $\sail$ containing a point $\vec x\in\sail$,
    note that $\varphi(\vec x)<\det S(F)$ and apply Theorem\,\ref{t:boundedness_of_F_and_SF}.
  \end{proof}

  \begin{lemma} \label{l:edges_are_short}
    Let the determinants of the facets of a sail $\sail$ be bounded by a constant $D$. Let
    $\vec v$ be a vertex of $\sail$ with $v_1=\ldots=v_n$ and let $\varphi(\vec v)\geq\mu>0$.
    Then, the (Euclidean) lengths of all the edges incident to the vertex $\vec v$ are bounded
    by a constant $D_{\textup{vert}}$ depending only on $n$, $D$ and $\mu$.
  \end{lemma}

  \begin{proof}
    Due to Theorem\,\ref{t:boundedness_of_F_and_SF}, there exists a constant $D'$
    depending only on $n$ and $D$, such that $\det S(F)\leq D'$ for each facet $F$ of the
    sail $\sail$.

    On the other hand, $v_1=\ldots=v_n\geq\mu^{1/n}$. Hence there exists a constant
    $D_{\textup{vert}}=D_{\textup{vert}}(n,D',\mu)$ such that if an edge incident to the
    vertex $\vec v$ has length larger than $D_{\textup{vert}}$, then for each facet $F$
    incident to this edge $\det S(F)>D'$.

    Therefore the lengths of all the edges incident to $\vec v$ should not exceed
    $D_{\textup{vert}}$.
  \end{proof}

  The following lemma is an obvious corollary of
  Statement\,\ref{st:polar_correspondence_Klein_polyhedra} and Definitions\,\ref{d:detF} and
  \ref{d:detStv}.

  \begin{lemma} \label{l:faces_and_vertices_of_K_ast_for_K_with_bounded_dets}
    If the facets and the edge stars of the vertices of a sail $\sail$ have determinants
    bounded by a constant $D$, then there exists a constant $D'$ depending only on $n$ and
    $D$ such that \\
    (a) each face of $K^\circ$ has not more than $D'$ vertices; \\
    (b) the number of facets of $K^\circ$ incident to a vertex of $K^\circ$ is bounded by
    $D'$.
  \end{lemma}

  \section{Boundedness away from zero of the form $\varphi(\vec x)$ in the positive orthant}

  As before, we suppose that the lattice $\La$ is irrational.

  \begin{lemma} \label{l:positiveness_in_the_positive_orthant}
    If the facets and the edge stars of the vertices of a sail $\sail$ have determinants
    bounded by a constant $D$, then there exists a constant $\mu>0$ depending only on $n$
    and $D$ for which
    \[ \inf_{\vec v}(\varphi(\vec v))\geq\mu, \]
    where the infimum is taken over all vertices of the sail $\sail$.
  \end{lemma}

  \begin{proof}
    It is easy to show that, if the boundary of the positive orthant contains nonzero points
    of the lattice $\La^\ast$, then the sail $\sail$ has an unbounded facet (see, e.g.,
    \cite{german_norm_minima_I}). But all facets of $\sail$ have bounded determinants,
    hence, there are no such points. Thus, by Statement\,\ref{st:Kast_is_a_subset_of_Kcirc},
    $K^\ast\subseteq K^\circ$. On the other hand, the integer distances from facets of $K$ to
    $\vec 0$ do not exceed $D$, hence, by
    Lemma\,\ref{l:if_distance_is_D_then_the_dual_lattice_gets_divided_by_D}, all vertices of
    $K^\circ$ lie in the lattice $(D!)^{-1}\La^\ast$. Applying the Corollary of
    Statement\,\ref{st:polar_correspondence_Klein_polyhedra}, we get that
    $D!\cdot K^\circ\subseteq K^\ast\subseteq K^\circ$.

    In virtue of Lemma\,\ref{l:boundedness_of_detFv}, the determinants of facets of $K^\circ$
    are bounded by $D^{n-1}$ and, thus, the determinants of facets of $D!\cdot K^\circ$
    are bounded by $D^{n-1}(D!)^n$. Let us prove the existence of a constant $D'$ depending
    only on $n$ and $D$ that bounds the determinants of facets of $K^\ast$. Due to the
    inclusion $D!\cdot K^\circ\subseteq K^\ast$, it suffices to show that the number of the
    facets of $K^\circ$ cut off by an arbitrary supporting hyperplane of $D!\cdot K^\circ$
    (including those that have nonempty intersection with this hyperplane) is bounded by a
    constant, which depends only on $n$ and $D$. Moreover, due to
    Lemma\,\ref{l:faces_and_vertices_of_K_ast_for_K_with_bounded_dets}, it suffices to
    consider only hyperplanes that are the affine hulls of facets of $D!\cdot K^\circ$. Let
    $F$ be a facet of $D!\cdot K^\circ$ and let $\aff(F)$ denote its affine hull. Obviously,
    the plane $\aff(F)$ contains an $(n-1)$--dimensional sublattice of $\La^\ast$, hence,
    the lattice $\La^\ast$, as well as the lattice $(D!)^{-1}\La^\ast$, can be split into
    $(n-1)$--dimensional layers parallel to $\aff(F)$. Consider now a facet $G$ of $K^\circ$
    such that the combinatorial distance between $(D!)^{-1}F$ and $G$ equals $k$ (here we
    call two different facets neighboring and we define the combinatorial distance between
    them to equal $1$, if they have at least one common point). It follows from the fact that
    all vertices of $K^\circ$ belong to $(D!)^{-1}\La^\ast$ that there are at least $k-2$
    layers of the lattice $(D!)^{-1}\La^\ast$ parallel to $\aff(F)$ such that their affine
    hulls strictly separate the facet $G$ from the facet $(D!)^{-1}F$. But since
    $\det F\leq D^{n-1}(D!)^n$ and $\det\La^\ast=1$, the number of layers of the lattice
    $(D!)^{-1}\La^\ast$ between $(D!)^{-1}\aff(F)$ and $\aff(F)$ is less than
    $D^{n-1}(D!)^{n+1}$. Therefore, applying
    Lemma\,\ref{l:faces_and_vertices_of_K_ast_for_K_with_bounded_dets}, we get that the
    number of facets of $K^\circ$ cut off by $\aff(F)$ is indeed bounded by a constant, which
    depends only on $n$ and $D$.

    Thus, all the facets of $K^\ast$ have determinants bounded by a constant $D'$ depending
    only on $n$ and $D$. Hence, by Theorem\,\ref{t:boundedness_of_F_and_SF}, there exists a
    constant $D''$, which also depends only on $n$ and $D$, such that
    $\varphi(\vec v)\geq(D'')^{-1}$ for each vertex $\vec v$ of $(K^\ast)^\ast=K$. It remains
    to set $\mu=(D'')^{-1}$.
  \end{proof}

  \section{The logarithmic plane}

  Let us denote the positive orthant by $\mathcal O_+$.
  Consider the two mappings:
  \[ \begin{array}{l}
       \pi_1:\mathcal O_+\rightarrow\big\{\vec x\in\mathcal O_+\,\big|\ \varphi(\vec x)=1\big\},
       \vphantom{\Big|} \\ \vphantom{\Big|}
       \pi_1(\vec x)=\vec x\cdot(\varphi(\vec x))^{-1/n}
     \end{array} \]
  and
  \[ \begin{array}{l}
       \pi_2:\big\{\vec x\in\mathcal O_+\,\big|\ \varphi(\vec x)=1\big\}\rightarrow\R^{n-1},
       \vphantom{\Big|} \\ \vphantom{\Big|}
       \pi_2(\vec x)=\left( \ln(x_1),\ldots,\ln(x_{n-1}) \right)
     \end{array} \]
  and their composition
  \[ \pilog=\pi_2\circ\pi_1. \]

  The image of a sail $\sail$ under the mapping $\pilog$ generates a partitioning of $\R^{n-1}$
  into cells being curvilinear polyhedra. Each cell is the image of some facet of the
  sail. Accordingly, each vertex of the partitioning is the image of some vertex of the sail.

  We denote this partitioning by $\partition$.

  \begin{lemma} \label{l:R_r_system}
    Suppose that the determinants of the facets of a sail $\sail$ are bounded by $D$ and
    there exists a constant $\mu>0$ such that, for each vertex $\vec v$ of the sail $\sail$,
    \[ \varphi(\vec v)\geq\mu. \]
    Then there exists a constant $D'$ depending only on $n$, $D$ and $\mu$, such that each
    ball $B\subset\pilog(\sail)$ of radius $D'$ contains a cell of the partitioning $\partition$.
  \end{lemma}

  \begin{proof}
    Consider an arbitrary vertex $\vec v$ of the sail $\sail$. Applying an appropriate
    hyperbolic rotation we can assume without loss of generality that
    $v_1=\ldots=v_n\geq\mu^{1/n}$. Then by Lemma\,\ref{l:edges_are_short}, the lengths of all
    the edges incident to $\vec v$ are bounded by some constant
    $D_{\textup{vert}}=D_{\textup{vert}}(n,D,\mu)$. Besides each facet of the sail has not
    more than $D_1=D_1(n,D)$ vertices. Therefore there exists such a constant
    $D_2=D_2(n,D,\mu)$, that all the facets incident to the vertex $\vec v$ are contained in
    a cube of sidelength $D_2$. But the values of the form $\varphi$ in all the points of the
    sail are not less than $\mu$. Hence, there exists a constant $D_3=D_3(n,D,\mu)$ such that
    the images under the mapping $\pilog$ of all the facets incident to the vertex $\vec v$
    are contained in a ball of radius $D_3$. Thus, for each facet $F$ of the sail $\sail$ the
    cell $\pilog(F)$ is contained in a ball of radius $D_3$.

    If we now consider an arbitrary ball $B\subset\pilog(\sail)$ of radius $D'=2D_3$, then
    the cell containing the center of $B$ is contained in a ball of radius $D_3$, which, in
    its turn, is contained in $B$.
  \end{proof}

  \section{Proof of Theorem\,\ref{t:the_main_one}}

  As before, we denote by $S(F)$ the intersection of the affine hull of a sail facet $F$
  and the positive orthant.

  {\bf 1.}
  The implication $(1)\,\&\,(2)\implies(3)$ has a rather simple proof. Consider an
  arbitrary vertex $\vec v$ of the sail $\sail$. Applying an appropriate hyperbolic
  rotation we can assume without loss of generality that $v_1=\ldots=v_n\geq\mu^{1/n}$.
  Then, by Lemma\,\ref{l:edges_are_short}, the lengths of all edges incident to $\vec v$ are
  bounded by a constant $D_{\textup{vert}}=D_{\textup{vert}}(n,\mu)$. This implies that the
  number of edges incident to $\vec v$ does not exceed some constant
  $m_{\textup{vert}}=m_{\textup{vert}}(n,\mu)$, because otherwise a ball of radius
  $D_{\textup{vert}}$ contains too many lattice points. Now, it is obvious that
  $\det\starv\leq D'=D'(D_{\textup{vert}},m_{\textup{vert}})=D'(n,\mu)$.

  {\bf 2.}
  The proof of the implication $(3)\implies (2)$ is a bit more difficult. We assume that the
  facets and the edge stars of the vertices of the sail $\sail$ have determinants bounded by
  a constant $D$. By Lemma\,\ref{l:positiveness_in_the_positive_orthant}, there exists a
  constant $\mu=\mu(n,D)>0$ such that, for each vertex $\vec v$ of the sail $\sail$, we have
  $ \varphi(\vec v)\geq\mu$, i.e. the conditions of Lemma\,\ref{l:R_r_system} are satisfied.

  Consider an arbitrary orthant $\mathcal O$ different from $\mathcal O_+$ and
  $-\mathcal O_+$ and an arbitrary facet $F$ of the sail corresponding to this orthant.
  Applying an appropriate hyperbolic rotation, we can assume that the facet $F$ is orthogonal
  to the bisector line of the orthant $\mathcal O$. We set
  \[ Q(T)=\big\{ \vec x\in\R^n\, \big|\ \max(|x_1|,\ldots,|x_n|)<T \big\}, \]
  \[ Q_+(T)=\big\{ \vec x\in Q(T)\, \big|\ x_i>0,\ i=1,\ldots,n \big\} \]
  and
  \[ T_0=n^{-1/2}(\det F)^{1/n}. \]
  Clearly, $Q(T_0)\cap\mathcal O\cap\La=\{\vec 0\}$.

  By virtue of Lemmas \ref{l:phi_is_bounded_on_sail} and \ref{l:R_r_system}, there exists a
  constant $T_1=T_1(n,D)$ such that the set $\pilog(\sail\cap Q_+(\sqrt{T_1}\,))$ contains a
  cell of the partitioning $\partition$ and, hence, a vertex of this partitioning. This means
  that $Q_+(\sqrt{T}\,)$ contains a vertex $\vec v$ of the sail $\sail$ for any $T\geq T_1$.
  Consider the parallelepiped
  \[ P(\vec v,T)=Q_+(T)\cap(\vec v+\mathcal O) \]
  for $T\geq T_1$. Lemmas \ref{l:phi_is_bounded_on_sail}, and \ref{l:R_r_system} imply the
  existence of a constant $T_2\geq T_1$, which also depends only on $n$ and $D$, such that
  the set $\pilog(\sail\cap P(\vec v,T_2))$ contains a cell and, hence, a vertex of the
  partitioning $\partition$. This means that $P(\vec v,T)$ contains a vertex of the sail
  $\sail$ different from $\vec v$ for any $T\geq T_2$.

  But the parallelepiped $P(\vec v,T_0)-\vec v$ is contained in the parallelepiped
  $Q(T_0)\cap\mathcal O$ and $Q(T_0)\cap\mathcal O\cap\La=\{\vec 0\}$. Hence, $T_0<T_2$,
  which means that $\det F$ is bounded by a constant depending only on $n$ and $D$.
  \qed


\begin{thebibliography}{99}
    \bibitem{german_norm_minima_II} \textsc{O. N. German},
        \textit{Klein polyhedra and norm minima of lattices}.
        Doklady Mathematics {\bf 406}:3 (2006), 38--41.
    \bibitem{erdos_gruber_hammer} \textsc{P. Erd\"{o}s, P. Gruber, J. Hammer},
        \textit{Lattice Points}.
        Pitman Monographs and Surveys in Pure and Applied
        Mathematics, {\bf 39}. Longman Scientific \& Technical, Harlow (1989).
    \bibitem{klein} \textsc{F. Klein},
        \textit{Uber eine geometrische Auffassung der gewohnlichen Kettenbruchentwichlung}.
        Nachr. Ges. Wiss. Gottingen, {\bf 3} (1895), 357--359.
    \bibitem{german_norm_minima_I} \textsc{O. N. German},
        \textit{Sails and norm minima of lattices}.
        Mat. Sb. {\bf 196}:3 (2005), 31--60;
        English transl., Russian Acad. Sci. Sb. Math. {\bf 196}:3 (2005), 337--367.
    \bibitem{moussafir_A_polyhedra} \textsc{J.--O. Moussafir},
        \textit{Convex hulls of integral points}.
        Zapiski nauch. sem. POMI, {\bf 256} (2000).
    \bibitem{arnold_mccme} \textsc{V. I. Arnold},
        \textit{Continued fractions}.
        Moscow: Moscow Center of Continuous Mathematical Education (2002).
    \bibitem{arnold_preface} \textsc{V. I. Arnold},
        \textit{Preface}.
        Amer. Math. Soc. Transl., {\bf 197}:2 (1999), ix--xii.
    \bibitem{korkina_2dim} \textsc{E. I. Korkina},
        \textit{Two--dimensional continued fractions. The simplest examples}.
        Proc. Steklov Math. Inst. RAS, {\bf 209} (1995), 143--166.
    \bibitem{bonnesen_fenchel} \textsc{T. Bonnesen, W. Fenchel},
        \textit{Theorie der konvexen K\"orper}.
        Berlin: Springer (1934).
    \bibitem{grunbaum} \textsc{B. Gr\"unbaum},
        \textit{Convex polytopes}.
        London, New York, Sydney: Interscience Publ. (1967).
    \bibitem{mcmullen_shephard} \textsc{P. McMullen, G. C. Shephard},
        \textit{Convex polytopes and the upper bound conjecture}.
        Cambridge (GB): Cambridge University Press (1971).
    \bibitem{ewald} \textsc{G. Ewald},
        \textit{Combinatorial convexity and algebraic geometry}.
        Sringer--Verlag New York, Inc. (1996).
    \bibitem{borevich_shafarevich} \textsc{Z. I. Borevich, I. R. Shafarevich},
        \textit{Number theory}.
        NY Academic Press (1966).
    \bibitem{cassels_swinnerton_dyer} \textsc{J. W. S. Cassels, H. P. F. Swinnerton--Dyer},
        \textit{On the product of three homogeneous linear forms and indefinite ternary
        quadratic forms}.
        Phil. Trans. Royal Soc. London, {\bf A 248} (1955), 73--96.
    \bibitem{skubenko_3dim} \textsc{B. F. Skubenko},
        \textit{Minima of a decomposable cubic form of three variables}.
        Zapiski nauch. sem. LOMI, {\bf 168} (1988).
    \bibitem{skubenko_ndim} \textsc{B. F. Skubenko},
        \textit{Minima of decomposable forms of degree $n$ of $n$ variables for $n\geq 3$}.
        Zapiski nauch. sem. LOMI, {\bf 183} (1990).
    \bibitem{lachaud_unseen} \textsc{G. Lachaud},
        \textit{Voiles et Poly\`{e}dres de Klein}.
        Act. Sci. Ind., Hermann (2002).
    \bibitem{danzer_grunbaum_klee} \textsc{L. Danzer, B. Gr\"unbaum, V. Klee},
        \textit{Helly's Theorem and its relatives}.
        in Convexity (Proc. Symp. Pure Math. 7) 101--180,
        AMS, Providence, Rhode Island (1963).
  \end{thebibliography}
\end{document}